\documentclass[12pt]{article}
\usepackage{graphicx}
\usepackage{amsmath,amsthm,amssymb,enumerate}
\usepackage{euscript,mathrsfs}
\usepackage[left=2cm,right=2cm,top=3.5cm,bottom=3.5cm]{geometry}
\usepackage{color}
\catcode`\@=11 \@addtoreset{equation}{section}

\catcode`\@=12

\allowdisplaybreaks

\newtheorem{Theorem}{Theorem}[section]
\newtheorem{Proposition}[Theorem]{Proposition}
\newtheorem{Lemma}[Theorem]{Lemma}
\newtheorem{Corollary}[Theorem]{Corollary}

\theoremstyle{definition}
\newtheorem{Definition}[Theorem]{Definition}

\newtheorem{Remark}[Theorem]{Remark}

\newcommand{\bTheorem}[1]{
\begin{Theorem} \label{T#1} }
\newcommand{\eT}{\end{Theorem}}

\newcommand{\bProposition}[1]{
\begin{Proposition} \label{P#1}}
\newcommand{\eP}{\end{Proposition}}

\newcommand{\bLemma}[1]{
\begin{Lemma} \label{L#1} }
\newcommand{\eL}{\end{Lemma}}

\newcommand{\bCorollary}[1]{
\begin{Corollary} \label{C#1} }
\newcommand{\eC}{\end{Corollary}}

\newcommand{\bRemark}[1]{
\begin{Remark} \label{R#1} }
\newcommand{\eR}{\end{Remark}}

\newcommand{\bDefinition}[1]{
\begin{Definition} \label{D#1} }
\newcommand{\eD}{\end{Definition}}

\newcommand{\avg}[1]{ \left< #1 \right>_\Gamma }
\newcommand{\vcg}[1]{{\pmb #1}}

\newcommand{\avo}[1]{ \left< #1 \right> }

\newcommand{\bFormula}[1]{
\begin{equation} \label{#1}}
\newcommand{\eF}{\end{equation}}

\newcommand{\TS}{\Delta t}

\newcommand{\Divh}{{\rm div}_h}
\newcommand{\Gradh}{\nabla_h}

\newcommand{\Ghi}{\Gamma_{h, {\rm int}}}

\newcommand{\Ov}[1]{\overline{#1}}
\newcommand{\av}[1]{ \left\{ #1 \right\}}

\newcommand{\DC}{C^\infty_c}

\newcommand{\aleq}{\stackrel{<}{\sim}}
\newcommand{\ageq}{\stackrel{>}{\sim}}

\newcommand{\vr}{\varrho}

\newcommand{\vu}{\vc{u}}
\newcommand{\vc}[1]{{\bf #1}}

\newcommand{\Div}{{\rm div}_x}
\newcommand{\Grad}{\nabla_x}

\newcommand{\dx}{\,{\rm d} {x}}

\newcommand{\dt}{\,{\rm d} t }

\newcommand{\ju}[1]{\left[ \left[ #1 \right] \right]}

\newcommand{\intO}[1]{\int_{\Omega} #1 \ \dx}
\newcommand{\intOh}[1]{\int_{\Omega_h} #1 \ \dx}

\newcommand{\intE}[1]{\int_E #1 \ \dx}
\newcommand{\intG}[1]{\int_{\Gamma} #1 \ {\rm dS}_x}

\newcommand{\cred}{\color{black}}

\definecolor{Cgrey}{rgb}{0.85,0.85,0.85}
\definecolor{Cblue}{rgb}{0.50,0.85,0.85}
\definecolor{Cred}{rgb}{1,0,0}
\definecolor{fancy}{rgb}{0.10,0.85,0.10}

\newcommand\Cbox[2]{%
    \newbox\contentbox%
    \newbox\bkgdbox%
    \setbox\contentbox\hbox to \hsize{%
        \vtop{
            \kern\columnsep
            \hbox to \hsize{%
                \kern\columnsep%
                \advance\hsize by -2\columnsep%
                \setlength{\textwidth}{\hsize}%
                \vbox{
                    \parskip=\baselineskip
                    \parindent=0bp
                    #2
                }%
                \kern\columnsep%
            }%
            \kern\columnsep%
        }%
    }%
    \setbox\bkgdbox\vbox{
        \color{#1}
        \hrule width  \wd\contentbox %
               height \ht\contentbox %
               depth  \dp\contentbox
        \color{black}
    }%
    \wd\bkgdbox=0bp%
    \vbox{\hbox to \hsize{\box\bkgdbox\box\contentbox}}%
    \vskip\baselineskip%
}

\date{}

\begin{document}


\title{Convergence of a mixed finite element--finite volume scheme for the isentropic Navier-Stokes system
via dissipative measure-valued solutions}

\author{Eduard Feireisl
\thanks{The research of E.F.~leading to these results has received funding from the
European Research Council under the European Union's Seventh
Framework Programme (FP7/2007-2013)/ ERC Grant Agreement
320078. The Institute of Mathematics of the Academy of Sciences of
the Czech Republic is supported by RVO:67985840.}\and M\'aria Luk{\' a}{\v c}ov{\' a}-Medvid'ov{\' a}
\thanks{The research of M.L.-M. has been supported by the German Science Foundation under the
grants LU 1470/2-3 and the Collaborative Research Centers TRR 146 and TRR 165.}
}

\date{\today}

\maketitle

\bigskip

\centerline{Institute of Mathematics of the Academy of Sciences of the Czech Republic}

\centerline{\v Zitn\' a 25, CZ-115 67 Praha 1, Czech Republic}
\medskip
\centerline{Institute of Mathematics, Johannes Gutenberg-University Mainz}

\centerline{Staudingerweg 6, 55 099 Mainz, Germany}

\bigskip

\begin{abstract}

We study convergence
of a mixed finite element--finite volume numerical scheme for the isentropic Navier-Stokes system under the full range of the adiabatic exponent.
We establish suitable stability and consistency estimates and show that the Young measure generated by numerical solutions
represents a dissipative measure-valued solutions of the limit system. {In particular},
using the recently established weak--strong uniqueness principle in the class of dissipative measure-valued solutions we show that the numerical solutions converge
strongly to a strong solutions of the limit system as long as the latter exists.

\end{abstract}


\newpage
\tableofcontents

\section{Introduction}
\label{i}

Time evolution of the density $\vr = \vr(t,x)$ and the velocity $\vu = \vu(t,x)$ of a compressible barotropic viscous fluid can be described
by the Navier--Stokes system
\begin{eqnarray}
\label{i1}
\partial_t \vr  + \Div (\vr \vu) &=& 0,
\\ \label{i2}
\partial_t (\vr \vu) + \Div (\vr \vu \otimes \vu) + \Grad p(\vr) &=& \Div \mathbb{S} (\Grad \vu),\\
\label{i3}
\mathbb{S} (\Grad \vu) &=& \mu \left( \Grad \vu + \Grad^t \vu - \frac{2}{3} \Div \vu \mathbb{I} \right) + \eta \Div \vu \mathbb{I}.
\end{eqnarray}
We assume the fluid is confined to a bounded physical domain $\Omega \subset R^3$, where the velocity satisfies the no-slip boundary conditions
\begin{equation} \label{i4}
\vu|_{\partial \Omega} =0.
\end{equation}
For the sake of simplicity, we ignore the effect of external forces in the momentum equation (\ref{i2}).

{\cred
In the literature there is a large variety of efficient numerical methods developed for the compressible Euler and Navier-Stokes equations. The most classical of them are the finite volume methods, see, e.g., \cite{feist1}, \cite{kroener}, \cite{tadmor-ns}, the methods based on a suitable
combination of the finite volume and finite element methods \cite{dolejsi}, \cite{ffl}, \cite{fflw}, \cite{gal1}, \cite{gal2},
or the discontinous Galerkin schemes,  e.g.~\cite{feist2}, \cite{feist3} and the references therein.
Although these methods are frequently used for many physical or engineering applications, there are only partial theoretical results available
concerning their analysis for the compressible Euler or
Navier-Stokes systems. We refer to the works of Tadmor et al.~\cite{tadmor2}, \cite{tadmor5}, \cite{tadmor} for entropy stability in the context of hyperbolic balance laws and to the works of Gallou\"et et al.~\cite{gal1}, \cite{gal2} for the stability analysis of the
mixed finite volume--finite element methods based on the Crouzeix-Raviart elements for compressible viscous flows.
In \cite{rohde} Jovanovi\'{c} and Rohde obtained the error estimate for entropy dissipative finite volume methods applied to nonlinear hyperbolic balance laws under (a rather restrictive) assumption of the global existence of a bounded, smooth exact solution.}

Our goal in this paper is to study convergence of solutions to the numerical scheme proposed originally by Karlsen and Karper \cite{KarKar3},
\cite{KarKar2}, \cite{KarKar1}, \cite{Karp}
to solve problem (\ref{i1}--\ref{i4}) in polygonal (numerical) domains, and later modified in
\cite{FeKaMi} to accommodate approximations of smooth physical domains. The scheme is implicit and of mixed type, where the convective terms are approximated via
upwind operators, while the viscous stress is handled by means of the Crouzeix--Raviart finite element method. As shown by Karper \cite{Karp} and in
\cite{FeKaMi},
the scheme provides a family of numerical solutions containing a sequence that converges to a weak solution of the Navier-Stokes system as the discretization
parameters tend to zero. Recently, Gallou{\"e}t et al. \cite{GalHerMalNov} established rigorous error estimates on condition that the limit problem admits
a smooth solution. Numerical experiments illustrating theoretical predictions have been performed in \cite{FLNNS}.

We consider the problem under physically realistic assumptions, where theoretical results are still in short supply. In particular, our results cover
completely the \emph{isentropic} pressure--density state equation
\begin{equation} \label{i5}
p(\vr) = a \vr^{\gamma}, \ 1 < \gamma < 2.
\end{equation}
Note that the assumption $\gamma < 2$ is not restrictive in this context as the largest physically relevant exponent is $\gamma = \frac{5}{3}$.
Let us remark that the available theoretical results concerning global-in-time existence of \emph{weak} solutions cover only the case
$\gamma > \frac{3}{2}$ \cite{FNP}, see also the recent result by Plotnikov and Weigant  \cite{PloWei} for the borderline case in the 2D setting.
{Similarly}, the
error estimates obtained by Gallou{\"e}t et al. \cite{GalHerMalNov} provide convergence under the same conditions yielding explicit convergence rates for
$\gamma > \frac{3}{2}$ and mere boundedness of the numerical solutions in the limit case $\gamma = \frac{3}{2}$.

Our goal is to establish convergence of the numerical solutions in the full range of the adiabatic exponent $\gamma$ specified in (\ref{i5}).
The main idea is to use the concept of \emph{dissipative measure-valued solution} to problem (\ref{i1}--\ref{i4}) introduced recently in {\cite{FGSWW1}}, \cite{GSWW}.
These are, roughly speaking, measure-valued solutions satisfying, in addition, an energy inequality in which the dissipation defect
measure dominates the concentration remainder in the equations. Although very general, a dissipative measure-valued solution coincides
with the strong solution of the same initial-value problem as long as the latter exists, see \cite{FGSWW1}. Our approach is based on the following steps:

\begin{itemize}
\item
We recall the numerical energy balance identified in Karper's original paper.
\item
We use the energy estimates to show stability of the numerical method.
\item
A consistency formulation of the problem is derived involving numerical solutions and error terms vanishing with the time step
$\TS$ and the spatial
discretization parameter $h$ approaching zero.
\item
We show that the family of numerical solutions generates a dissipative measure-valued solution of the problem. Such a result is, of course,
of independent interest. {\cred As claimed recently by Fjordholm et al.~\cite{tadmor3}, \cite{tadmor4}  the dissipative measure-valued solutions yield,
{at least in the context of \emph{hyperbolic} conservation laws}, a more appropriate solution concept than the weak entropy solutions.}
\item
Finally, using the weak--strong uniqueness principle established in \cite{FGSWW1}, we infer that the numerical solutions converge (a.a.) pointwise
to the smooth solution of the limit problem as long as the latter exists.

\end{itemize}

The paper is organized as follows. The numerical scheme is introduced in Section \ref{N}. In Section \ref{S}, we recall the numerical counterpart of the
energy balance and derive stability estimates. In Section \ref{C}, we introduce a consistency formulation of the problem and estimate the numerical
errors. Finally, we show that the numerical scheme generates a dissipative measure-valued solution to the compressible Navier--Stokes system
and state our main convergence results in Section \ref{M}.

\section{Numerical scheme}
\label{N}

To begin, we introduce the notation necessary to formulate our numerical method.

\subsection{Spatial domain, mesh}

We suppose that $\Omega \subset R^3$ is a bounded domain. We consider a polyhedral approximation $\Omega_h$, where $\Omega_h$ is a polygonal domain,
\[
\Ov{\Omega}_h = \cup_{E^j \in E_h} E^j, \ {\rm int}[ E^i ] \cap {\rm int}[E^j] = \emptyset \ \mbox{for}\ i \ne j,
\]
where each $E^j \in E_h$ is a closed tetrahedron that can be obtained via the affine transformation
\[
E^j = h \mathbb{A}_{E^j} \tilde E + \vc{a}_{E^j}, \ \mathbb{A}_{E^j} \in R^{3 \times 3}, \ \vc{a}_{E^j} \in R^3,
\]
where $\tilde E$ is the reference element
\[
\tilde{E} = {\rm co} \left\{ [0,0,0], [1,0,0], [0,1,0], [0,0,1] \right\},
\]
and
where all eigenvalues of the matrix $\mathbb{A}_{E^j}$ are bounded above and below away from zero uniformly for $h \to 0$.
The family $E_h$ of all tetrahedra covering $\Omega_h$ is called \emph{mesh}, the positive number $h$ is the parameter of spatial discretization.
We write
\[
\begin{split}
a &\aleq b \Leftrightarrow a \leq c b, \ c > 0 \ \mbox{independent of}\ h, \\
a &\ageq b \Leftrightarrow a \geq c b, \ c > 0 \ \mbox{independent of}\ h, \\
a&=b \Leftrightarrow a \aleq b \ \mbox{and}\ a \ageq b.
\end{split}
\]

Furthermore, we suppose that:

\begin{itemize}
\item
a non-empty intersection of two elements $E^j$, $E^i$ is their common face, edge, or vertex;
\item
for all compact sets $K_i \subset \Omega$, $K_e \subset R^3 \setminus \Ov{\Omega}$ there is $h_0 > 0$ such that
\[
K_i \subset \Omega_h, \ K_e \subset R^3 \setminus \Ov{\Omega}_h \ \mbox{for all}\ 0 < h < h_0.
\]
\end{itemize}

The symbol $\Gamma_h$ denotes the set of all faces in the mesh. We distinguish exterior and interior faces:
\[
\Gamma_h = \Gamma_{h, {\rm int}} \cup \Gamma_{h, {\rm ext}},\
\Gamma_{h, {\rm ext}} = \left\{ \Gamma \in \Gamma_h \ \Big| \ \Gamma \subset \partial \Omega_h \right\},\
\Gamma_{h, {\rm int}} = \Gamma_h \setminus \Gamma_{h, {\rm ext}}.
\]

\subsection{Function spaces}

Our scheme utilizes spaces of piecewise smooth functions, for which we define the traces
\[
v^{\rm out} = \lim_{\delta \to 0} v(x + \delta \vc{n}_\Gamma),\
v^{\rm in} = \lim_{\delta \to 0} v(x - \delta \vc{n}_\Gamma), \ x \in \Gamma, \ \Gamma \in \Gamma_{h, {\rm int}},
\]
where $\vc{n}_\Gamma$ denotes the outer normal vector to the face $\Gamma \subset \partial E$. Analogously, we define
$v^{\rm in}$ for $\Gamma \subset \Gamma_{h,{\rm ext}}$. We simply write $v$ for $v^{\rm in}$ if no confusion arises. We also define
\[
\ju{v} = v^{\rm out} - v^{\rm in}, \ \avg{ v } = \frac{ v^{\rm out} + v^{\rm in} }{2}, \ \avg{v} = \frac{1}{|\Gamma|} \intG{v}.
\]

Next, we introduce the space of piecewise
constant functions
\[
Q_h (\Omega_h) = \left\{ v \in L^1(\Omega_h) \ \Big|\ v|_E = {\rm const} \in R \ \mbox{for any} \ E \in E_h \right\},
\]
with the associated projection
\[
\Pi^{Q}_h : L^1 (\Omega_h) \to Q_h (\Omega_h), \ \Pi^Q_h [v] = \left< v \right>_{E} = \frac{1}{|E|} \int_E v \ \dx, \ E \in E_h.
\]
We shall occasionally write
\[
\Pi^Q_h [v] = \avo{v}.
\]

Finally, we introduce the Crouzeix--Raviart finite element spaces
\[
V_h (\Omega_h) = \left\{ v \in L^2(\Omega_h) \ \Big| \ v|_E = \mbox{affine function} \ E \in E_h, \
\intG{ v^{\rm in} } = \intG{ v^{\rm out}} \ \mbox{for}\ \Gamma \in \Gamma_{h, {\rm int}} \right\},
\]
\[
V_{0,h} (\Omega_h) = \left\{ v \in V_h (\Omega_h) \ \Big| \ \intG{ v^{\rm in} } = 0 \ \mbox{for}\ \Gamma \in \Gamma_{h, {\rm ext}} \right\},
\]
along with the associated projection
\[
\Pi^V_h : W^{1,1}(\Omega_h) \to V_h (\Omega_h), \ \intG{ \Pi^V_h [v] } = \intG{v} \ \mbox{for any}\ \Gamma \in \Gamma_h.
\]
We denote by $\nabla_h v$, ${\rm div}_h v$ the piecewise constant functions resulting from the action of the corresponding differential
operator on $v$ on each fixed element in $E_h$,
\[
\nabla_h v \in Q_h(\Omega_h; R^3), \ \nabla_h v = \Grad v \ \mbox{for}\ E \in E_h,\
{\rm div}_h \vc{v} \in Q_h(\Omega_h), \ {\rm div}_h v = \Div v \ \mbox{for}\ E \in E_h.
\]

\subsection{Discrete time derivative, dissipative upwind}

For a given time step $\TS > 0$ and the (already known) value of the numerical solution $v^{k-1}_h$ at a given time level $t_{k-1} = (k-1) \TS$,
we introduce the discrete time derivative
\[
D_t v_h = \frac{ v^k_h - v^{k-1}_h }{\TS}
\]
to compute the numerical approximation $v^k_h$ at the level $t_k = t_{k-1} + \TS$.

To approximate the convective terms, we use the dissipative upwind operators introduced in \cite{FeKaMi} (see also \cite{FeiKaPok}), specifically,
\begin{equation} \label{N1}
\begin{split}
{\rm Up}[r_h, \vu_h] &= \underbrace{\av{ r_h } \avg{ \vu_h \cdot \vc{n}}}_{\rm convective \ part} - \frac{1}{2}
\underbrace{\max\{ h^\alpha; | \avg{ \vu_h \cdot \vc{n}} | \} \ju{ r_h }}_{\rm dissipative \ part}
\\ &= \underbrace{ r_h^{\rm out} [ \avg{\vu_h \cdot \vc{n}} ]^- +
r_h^{\rm in} [ \avg{\vu_h \cdot \vc{n}} ]^+}_{\rm standard \ upwind} - \frac{h^\alpha}{2}  \ju{r_h}  \chi \left( \frac{ \avg{\vu_h \cdot \vc{n} }}{h^\alpha} \right),
\end{split}
\end{equation}
where
\[
\chi(z) = \left\{ \begin{array}{l} 0 \ \mbox{for}\ z < -1, \\  z + 1 \ \mbox{if} \ -1 \leq z \leq 0, \\
1 - z \ \mbox{if} \ 0 < z \leq 1, \\
0 \ \mbox{for}\ z > 1.
\end{array} \right.
\]

\subsection{Numerical scheme}

Given the initial data
\begin{equation} \label{N2}
\vr^0_h \in Q_h (\Omega_h), \ \vu^0_h \in V_{0,h} (\Omega_h; R^3),
\end{equation}
and the numerical solution
\[
\vr^{k-1}_h \in Q_h (\Omega_h), \ \vu^{k-1}_h \in V_{0,h} (\Omega_h; R^3),\ k \geq 1,
\]
the value $[\vr^k_h, \vu^k_h] \in Q_h (\Omega_h) \times V_{0,h} (\Omega_h; R^3)$ is obtained as a solution of the following system of equations:

\begin{equation} \label{N3}
\intOh{ D_t \vr^k_h \phi } - \sum_{\Gamma \in \Ghi} \intG{ {\rm Up}[\vr^k_h, \vu^k_h] \ju{\phi} } = 0
\end{equation}
for any $\phi \in Q_h(\Omega_h)$;
\begin{equation} \label{N4}
\begin{split}
\intOh{ D_t \left( \vr^k_h \avo{ \vu^k_h } \right) \cdot \vcg{\phi} } &- \sum_{\Gamma \in \Ghi} \intG{ {\rm Up}[\vr^k_h \avo{ \vu^k_h} , \vu^k_h]
\cdot \ju{ \avo{ \vcg{\phi} } } }- \intOh{ p(\vr^k_h) \Divh \vcg{\phi} }\\
&+ \mu \intOh{ \Gradh \vu^k_h : \Gradh \vcg{\phi} } + \left( \frac{\mu}{3} + \eta \right) \intOh{ \Divh \vu^k_h \Divh \vcg{\phi} } =
0
\end{split}
\end{equation}
for any $\vcg{\phi} \in V_{0,h}(\Omega_h; R^3)$. The specific form of the viscous stress in (\ref{N4}) reflects the fact that the viscosity coefficients
are constant.

It was shown in {\cite{Karp} (see also \cite[Part II]{FeiKaPok})} that system (\ref{N3}), (\ref{N4}) is solvable for any choice of the initial data (\ref{N2}). In addition,
$\vr^k_h > 0$ whenever $\vr^0_h > 0$. In general, the solution $[\vr^k_h, \vu^k_h]$ may not be uniquely determined by $[\vr^{k-1}_h, \vu^{k-1}_h]$
unless the time step $\TS$ is conveniently adjusted by a CFL type condition. We make more comments on this option in Remark \ref{RC3} below.

As shown in {\cite{FeKaMi} (see also \cite[Part II]{FeiKaPok})}, the family of numerical solutions converges, up to a suitable subsequence, to a weak solution of the Navier-Stokes
system (\ref{i1}--\ref{i4}) as $h \to 0$ if
\begin{itemize}
\item
the time step is adjusted so that $\TS \approx h$;
\item
the viscosity coefficients satisfy $\mu > 0$, $\eta \geq 0$,
\item
the pressure satisfies
\[
p(\vr) = a \vr^{\gamma} + b \vr, \ a,b > 0, \ \gamma > 3.
\]
\end{itemize}
If the limit solution of the Navier--Stokes system is smooth, then qualitative error estimates can be derived on condition
that $p$ satisfies (\ref{i5}) with $\gamma \geq 3/2$, {see Gallou{\"e}t et al. \cite{GalHerMalNov}.
Unfortunately, many real world applications correspond to smaller adiabatic exponents, the most popular among them is the air with
$\gamma = 7/5$.
It is therefore of great interest to discuss convergence
of the scheme in the physically relevant range
$1 < \gamma < 2$.}

\section{Stability - energy estimates}
\label{S}

It is crucial for our analysis that the numerical scheme (\ref{N2}--\ref{N4}) admits a certain form of total energy balance.
For the pressure potential
\[
P(\vr) = \frac{a}{\gamma - 1} \vr^{\gamma}, \ P''(\vr) =  \frac{p'(\vr)}{\vr} = a \gamma \vr^{\gamma - 2},
\]
the \emph{total energy balance} reads
\begin{equation} \label{S1}
\begin{split}
&\intOh{ D_t \left[ \frac{1}{2} \vr^k_h |\avo{\vu^k_h}|^2 + P(\vr^k_h) \right] }
+ \intOh{  \left[\mu |\Gradh \vu^k_h |^2 + (\mu/3 + \eta)
|\Divh \vu^k_h |^2\right]  }
\\
&= - \frac{1}{2} \intOh{ P''(s^k_h) \frac{ \left( \vr^k_h - \vr^{k-1}_h \right)^2}{\TS} } - \intOh{\frac{\TS}{2} \vr^{k-1}_h
\left| \frac{ \avo{\vu^k_h} -  \avo{\vu^{k-1}_h} } {\TS} \right|^2  }
\\
&- \frac{h^\alpha}{2} \sum_{\Gamma \in \Ghi} \intG{
\ju{\vr^k_h} \ju{ P' (\vr^k_h )}\chi \left( \frac{ \avg{\vu^k_h \cdot \vc{n} }}{h^\alpha} \right) }
\\
&- \frac{1}{2} \sum_{\Gamma \in \Gamma_h} \intG{ P''(z^k_h) \ju{ \vr^k_h}^2 | \avg{ \vu^k_h \cdot \vc{n} } | }
\\
&- \frac{h^\alpha}{2} \sum_{\Gamma \in \Ghi} \intG{ \av{ \vr^k_h  } \cdot \ju{ \avo{\vu^k_h} }^2
\chi \left( \frac{ \avg{ \vu^k_h \cdot \vc{n} } }{h^\alpha} \right)  }
\\
&- \frac{1}{2} \sum_{\Gamma \in \Ghi } \intG{ \left( (\vr^k_h)^{\rm in} [ \avg{ \vu^k_h \cdot \vc{n} } ]^+ -
(\vr^k_h)^{\rm out} [ \avg{ \vu^k_h \cdot \vc{n}} ]^- \right) \ju{ \avo{\vu^k_h} }^2 },
\end{split}
\end{equation}
with
\[
s^k_h \in {\rm co}\{ \vr^k_h, \vr^{k-1}_h \}, \ z^k_h \in {\rm co}\{ (\vr^k)^{\rm in}, (\vr^k_h)^{\rm out} \},
\]
see \cite[Chapter 7, Section 7.5.4]{FeiKaPok}.
As the numerical densities are positive,
all terms on the right-hand side of (\ref{S1}) representing numerical dissipation
are non-positive. For completeness, we remark that the scheme conserves the total mass, specifically,
\begin{equation} \label{S1a}
\intOh{ \vr^k_h } = \intOh{ \vr^0_h }, \ k=1,2, \dots
\end{equation}

\subsection{Dissipative terms and the pressure growth}

It is easy to check that
\begin{equation} \label{S1b}
P''(z) (\vr_1 - \vr_2)^2 \geq a\gamma (\vr^{\gamma/2}_1 - \vr^{\gamma/2}_2 )^2 \ \mbox{whenever}\
z \in {\rm co} \{ \vr_1, \vr_2 \}, \ \vr_1, \vr_2 > 0, \ 1 < \gamma < 2.
\end{equation}
Indeed it is enough to assume $0 < \vr_1 \leq z \leq \vr_2$; whence
\[
P''(z) (\vr_1 - \vr_2)^2 \geq a \gamma \vr_2^{\gamma - 2} (\vr_1 - \vr_2)^2,
\]
and (\ref{S1b}) reduces to showing
\[
\vr_2^{\gamma/2 - 1} (\vr_2 - \vr_1) \geq (\vr^{\gamma/2}_2 - \vr^{\gamma/2}_1 )
\ \mbox{or, equivalently,}\ \vr_1 \vr_2^{\gamma/2 - 1} \leq \vr_1^{\gamma/2},
\]
where the last inequality follows immediately as $\vr_1 \leq \vr_2$, $1 < \gamma < 2$.

Consequently, the terms on the right-hand side of (\ref{S1}) representing the numerical dissipation and containing $P''$ satisfy
\begin{equation} \label{S2}
\begin{split}
\frac{1}{2} \intOh{ P''(s^k_h) \frac{ \left( \vr^k_h - \vr^{k-1}_h \right)^2}{\TS} }
&\geq \frac{a\gamma}{2} \intOh{  \frac{ \left( (\vr^k_h)^{\gamma/2} - (\vr^{k-1}_h)^{\gamma/2} \right)^2}{\TS} },\\
\frac{h^\alpha}{2} \sum_{\Gamma \in \Ghi} \intG{
\ju{\vr^k_h} \ju{ P' (\vr^k_h )}\chi \left( \frac{ \avg{\vu^k_h \cdot \vc{n} }}{h^\alpha} \right) } &\geq
\frac{a\gamma h^\alpha}{2} \sum_{\Gamma \in \Ghi} \intG{
\ju{ (\vr^k_h)^{\gamma / 2} }^2  \chi \left( \frac{ \avg{\vu^k_h \cdot \vc{n} }}{h^\alpha} \right) },\\
\frac{1}{2} \sum_{\Gamma \in \Gamma_h} \intG{ P''(z^k_h) \ju{ \vr^k_h}^2 | \avg{ \vu^k_h \cdot \vc{n} } | } &\geq
\frac{a\gamma}{2} \sum_{\Gamma \in \Gamma_h} \intG{ \ju{ (\vr^k_h )^{\gamma/2} }^2 | \avg{ \vu^k_h \cdot \vc{n} } | }.
\end{split}
\end{equation}
In particular, the energy balance (\ref{S1}) gives rise to
\begin{equation} \label{S3}
\begin{split}
&\intOh{ D_t \left[ \frac{1}{2} \vr^k_h |\avo{\vu^k_h}|^2 + P(\vr^k_h) \right] }
+ \intOh{  \left[\mu |\Gradh \vu^k_h |^2 + (\mu/3 + \eta)
|\Divh \vu^k_h |^2\right]  }
\\
& + a\intOh{  \frac{ \left( (\vr^k_h)^{\gamma/2} - (\vr^{k-1}_h)^{\gamma/2} \right)^2}{\TS} } + \TS \intOh{ \vr^{k-1}_h
\left| \frac{ \avo{\vu^k_h} -  \avo{\vu^{k-1}_h} } {\TS} \right|^2  }
\\
&+  a\sum_{\Gamma \in \Gamma_h} \intG{  \ju{ (\vr^k_h)^{\gamma/2} }^2 \max\left\{ h^\alpha ; | \avg{ \vu^k_h \cdot \vc{n} } | \right\} }
\\
&+ a{h^\alpha}\sum_{\Gamma \in \Ghi} \intG{ \av{ \vr^k_h  } \cdot \ju{ \avo{\vu^k_h} }^2
\chi \left( \frac{ \avg{ \vu^k_h \cdot \vc{n} } }{h^\alpha} \right)  }
\\
&+  \sum_{\Gamma \in \Ghi } \intG{ \left( (\vr^k_h)^{\rm in} [ \avg{ \vu^k_h \cdot \vc{n} } ]^+ -
(\vr^k_h)^{\rm out} [ \avg{ \vu^k_h \cdot \vc{n}} ]^- \right) \ju{ \avo{\vu^k_h} }^2 }  \aleq 0.
\end{split}
\end{equation}

\section{Consistency}
\label{C}

Our goal is to derive a consistency formulation for the discrete solutions satisfying (\ref{N3}), (\ref{N4}). To this end, it is convenient
to deal with quantities defined on $R \times \Omega_h$. Accordingly, we introduce
\begin{equation} \label{C1}
\vr_h (t, \cdot) = \vr^0_h \ \mbox{for} \ t < \TS, \ \vr_h(t, \cdot) = \vr^k_h \ \mbox{for}\ t \in [k \TS, (k+1) \TS),\ k = 1, 2, \dots,
\end{equation}
\begin{equation} \label{C2}
\vu_h (t, \cdot) = \vu^0_h \ \mbox{for} \ t < \TS, \ \vu_h(t, \cdot) = \vu^k_h \ \mbox{for}\ t \in [k \TS, (k+1) \TS), \ k = 1, 2, \dots,
\end{equation}
and
\begin{equation} \label{C3}
D_t v_h = \frac{ v(t, \cdot) - v(t - \TS, \cdot) }{\TS}, \ t > 0.
\end{equation}
For the sake of simplicity, we keep the time step $\TS$ constant, however, a similar ansatz obviously works also for $\TS = \TS_k$ adjusted
at each level of iteration.

A suitable consistency formulation of equation (\ref{N3}) reads
\begin{equation} \label{C10}
- \intOh{ \vr^0_h \varphi (0, \cdot) } = \int_0^T \intOh{ \left[  \vr_h \partial_t \varphi + \vr_h \vu_h \cdot \Grad \varphi \right] } \dt
+ \mathcal{O}(h^\beta), \ \beta > 0,
\end{equation}
for any test function $\varphi \in \DC([0, \infty) \times \Ov{\Omega}_h )$, where $\beta$ denotes a generic positive exponent, and,
accordingly, the remainder term $\mathcal{O}(h^\beta)$, that may depend also on the test function
$\varphi$, tends to zero as $h \to 0$. Similarly, we want to rewrite (\ref{N4}) in the form
\begin{equation} \label{C16}
\begin{split}
- \intOh{ \vr^0_h \avo{ \vu^0_h } \cdot  \vcg{\varphi}(0, \cdot) } &= \int_0^T \intOh{
\Big[ \vr_h \avo{\vu_h} \partial_t \vcg{\varphi} + \vr_h \avo{\vu_h} \otimes \vu_h : \Grad \vcg{\varphi} + p(\vr_h ) \Div \vcg{\varphi} \Big] } \ \dt\\
&- \int_0^T \intOh{ \Big[ \mu \nabla_h \vu_h : \Grad \vcg{\varphi} + (\mu/3 + \eta) {\rm div}_h \vu_h \cdot \Div \vcg{\varphi} \Big] } \dt
+ \mathcal{O}(h^\beta)
\end{split}
\end{equation}
for any $\vcg{\varphi} \in \DC([0, \infty) \times \Omega_h; R^3)$.

\subsection{Preliminaries, some useful estimates}

We collect certain well-known estimates used in the subsequent analysis. We refer to \cite[Part II, Chapters 8,9]{FeiKaPok} for the proofs.

\subsubsection{Discrete negative and trace estimates for piecewise smooth functions}

The following inverse inequality
\begin{equation} \label{S4}
\| v \|_{L^p(\Omega_h) } \aleq h^{ 3 \left( \frac{1}{p} - \frac{1}{q} \right) } \| v \|_{L^q(\Omega_h)}, \ 1 \leq q \leq p \leq \infty,
\end{equation}
holds for any $v \in Q_h(\Omega_h)$.

The trace estimates read
\begin{equation} \label{S4a}
\| v \|_{L^p(\Gamma)} \aleq h^{1/p} \| v \|_{L^p(E)} \ \mbox{whenever}\ \Gamma \subset \partial E, \ 1 \leq p \leq \infty
\end{equation}
for any $v \in Q_h(\Omega_h)$.

Finally, we report a discrete version of Poincar\' e's inequality
\begin{equation} \label{S4b}
\| v - \avo{v} \|_{L^2(E)} \equiv \| v - \Pi^Q_h [v] \|_{L^2(E)} \aleq h \| \nabla_h v \|_{L^2(E)}
\ \mbox{for any}\ v \in V_h(\Omega_h).
\end{equation}

\subsubsection{Sobolev estimates for broken norms}

We have
\begin{equation} \label{S5}
\| v \|_{L^6(\Omega_h)}^2 \aleq \sum_{\Ghi} \intG{ \frac{\ju{ v }^2}{h} } + \| v \|^2_{L^2(\Omega_h)}
\end{equation}
for any $v \in Q_h(\Omega_h)$.
In particular, we may combine the negative estimates (\ref{S4}) with (\ref{S5}) to obtain
\begin{equation} \label{S6}
\begin{split}
\| \vr_h \|_{L^\infty(\Omega_h)} &= \left( \left\| \vr_h^{\gamma/2} \right\|_{L^\infty(\Omega_h)} \right)^{2 / \gamma} \aleq
h^{-1/\gamma} \left( \left\| \vr_h^{\gamma/2} \right\|_{L^6(\Omega_h)}^2 \right)^{1 / \gamma}\\
&\aleq h^{-1/\gamma} \left( \sum_{\Ghi} \intG{ \frac{\ju{ \vr^{\gamma/2}  }^2}{h} }\right)^{1 / \gamma}
+ h^{-1/\gamma} \left(  \left\| \vr^{\gamma/2} \right\|_{L^2(\Omega_h)}^2  \right)^{1 / \gamma}\\
&\aleq h^{- \frac{2 + \alpha}{\gamma} } \left( \sum_{\Ghi} \intG{ h^\alpha {\ju{ \vr^{\gamma/2}  }^2} }\right)^{1 / \gamma}
+ h^{-1/\gamma}  \left\| \vr \right\|_{L^\gamma (\Omega_h)}
\end{split}
\end{equation}

Next, we have the discrete variant of Sobolev's inequality
\begin{equation} \label{S6b}
\| v \|^2_{L^6(\Omega_h)} \aleq \sum_{E \in E_h} \| \nabla_h v \|^2_{L^2(E; R^3)} \equiv \| \nabla_h v \|_{L^2(\Omega_h; R^3)}^2
\end{equation}
for any $v \in V_{0,h}(\Omega_h)$.

Finally, we recall the projection estimates for the Crouzeix--Raviart spaces
\begin{equation} \label{S7}
\left\| \Pi^V_h [v] - v \right\|_{L^q(\Omega_h)} + h \left\| \nabla_h \Pi^V_h [v] - \Grad v \right\|_{L^q(\Omega_h;R^3)}
\aleq h^j \| \nabla^j v \|_{L^q(\Omega_h; R^{3j})},\ j=1,2, \ 1 \leq q \leq \infty.
\end{equation}

\subsubsection{Upwind consistency formula}

We report the universal formula
\begin{equation} \label{C4}
\begin{split}
\intOh{ r \vu \cdot \Grad \phi } &= \sum_{\Gamma \in \Ghi} \intG{ {\rm Up}[r, \vu] \ju{F} }
\\
&+ \frac{h^\alpha}{2} \sum_{\Gamma \in \Ghi} \intG{ \ju{r} \ju{F} \chi \left( \frac{ \avg{\vu \cdot \vc{n} }}{h^\alpha} \right) }
\\
&+ \sum_{E \in E_h} \sum_{\Gamma_E \subset \partial E} \int_{\Gamma_E} (F - \phi) \ju{r} [ \avg{ \vu \cdot \vc{n} }]^- \ {\rm dS}_x
\\
&+ \sum_{E \in E_h} \sum_{\Gamma_E \subset \partial E} \int_{\Gamma_E} \phi r \Big(\vu \cdot \vc{n} - \avg{ \vu \cdot \vc{n} } \Big) \ {\rm dS}_x
+ \intOh{ r (F - \phi) \Divh \vu }
\end{split}
\end{equation}
for any
$r, F \in Q_h(\Omega_h)$, $\vu \in V_{0,h}(\Omega_h; R^3)$, $\phi \in C^1(\Omega_h)$, see \cite[Chapter 9, Lemma 7]{FeiKaPok}.

\subsection{Consistency formulation of the continuity method}

Our goal is to derive the consistency formulation (\ref{C10}) of the discrete equation of continuity (\ref{N3}).

\subsubsection{Time derivative}
\label{TD1}

We consider test functions of the form $\psi(t) \phi(x)$ to obtain
\[
\begin{split}
\int_0^T \intOh{ D_t (\vr_h) \avo{\psi \phi} } &\ \dt = \int_0^T \psi \intOh{ D_t (\vr_h) \phi } \ \dt\\
&= - \int_0^T \intOh{ \frac{ \psi(t + \TS) - \psi(t) }{\TS} \vr_h \phi } \ \dt - \frac{1}{\TS} \int_{-\TS}^0 \intOh{ \vr^0_h \psi (t + \TS) \phi }\ \dt
\end{split}
\]
whenever the function $\psi \in \DC[0, T)$ and $\TS$ is small enough so that the interval $[T - \TS, \infty)$ is not included in the support of $\psi$.
By means of the mean-value theorem we get that
\begin{equation} \label{C5}
\int_0^T \intOh{ D_t (\vr_h) \avo{\psi \phi} } \ \dt
= - \int_0^T \intOh{ \partial_t \psi \vr_h \phi } \ \dt  -  \intOh{ \vr^0_h \psi (0) \phi } + \mathcal{O}(h^\beta)
\end{equation}
for any $\phi \in C(\Omega_h)$, $\psi \in \DC[0, T)$. Note that the $\mathcal{O}(h)$ term depends on the second derivative of $\psi$.

\subsubsection{Convective term - upwind}

Relation (\ref{C4}) evaluated for $r = \vr^k_h$, $\vu = \vu^k_h$, $F = \avo{\phi}$, $\phi \in C^1(\Omega_h)$ gives rise to
\begin{equation} \label{C6}
\begin{split}
\intOh{ \vr^k_h \vu^k_h \cdot \Grad \phi } &= \sum_{\Gamma \in \Ghi} \intG{ {\rm Up}[\vr^k_h, \vu^k_h] \ju{\avo{\phi}} }
\\
&+ \frac{h^\alpha}{2} \sum_{\Gamma \in \Ghi} \intG{ \ju{\vr^k_h} \ju{\avo{\phi} } \chi \left( \frac{ \avg{\vu^k_h \cdot \vc{n} }}{h^\alpha} \right) }
\\
&+ \sum_{E \in E_h} \sum_{\Gamma_E \subset \partial E} \int_{\Gamma_E} (\avo{\phi} - \phi) \ju{\vr^k_h} [ \avg{ \vu^k_h \cdot \vc{n} }]^- \ {\rm dS}_x
\\
&+ \sum_{E \in E_h} \sum_{\Gamma_E \subset \partial E} \int_{\Gamma_E} \phi \vr^k_h \Big(\vu^k_h \cdot \vc{n} - \avg{ \vu^k_h \cdot \vc{n} } \Big) \ {\rm dS}_x
+ \intOh{ \vr_h (\avo{\phi} - \phi) \Divh \vu^k_h }.
\end{split}
\end{equation}

Using an elementary inequality
\begin{equation} \label{C6a}
\left| \vr_1 - \vr_2 \right| \leq \left| (\vr_1)^{\gamma/2} - (\vr_2)^{\gamma/2} \right| \left| (\vr_1)^{1 - \gamma/2} + (\vr_2)^{1 - \gamma/2} \right|,\
1 \leq \gamma \leq 2
\end{equation}
we get
\[
\begin{split}
\frac{h^\alpha}{2} &\left| \sum_{\Gamma \in \Ghi} \intG{ \ju{\vr^k_h} \ju{\avo{\phi} } \chi \left( \frac{ \avg{\vu^k_h \cdot \vc{n} }}{h^\alpha} \right) }
\right|
\aleq h^{1 + \alpha} \| \phi \|_{C^1(\Ov{\Omega}_h)} \left| \sum_{\Gamma \in \Ghi} \intG{ \ju{\vr^k_h}  } \right| \\
&\aleq
 h^{1 + \alpha} \| \phi \|_{C^1(\Ov{\Omega}_h)} \left( \sum_{\Gamma \in \Ghi} \intG{ \ju{(\vr^k_h)^{\gamma/2}}^2  } +
\sum_{\Gamma \in \Ghi} \intG{ \av{ (\vr^k_h)^{1 - \gamma/2} }^2  } \right),
\end{split}
\]
where, by virtue of (\ref{S3}),
\[
 h^{1 + \alpha} \| \phi \|_{C^1(\Ov{\Omega}_h)}  \sum_{\Gamma \in \Ghi} \intG{ \ju{(\vr^k_h)^{\gamma/2}}^2  } \leq
c(\phi) h g_k,\ \TS \sum_{k} g_k < \infty,
\]
and, in accordance with  (\ref{S1a}) and the trace estimates (\ref{S4a}),
\[
 h^{1 + \alpha} \| \phi \|_{C^1(\Ov{\Omega}_h)} \sum_{\Gamma \in \Ghi} \intG{ \av{ (\vr^k_h) ^{1 - \gamma/2} }^2  }
 \aleq h^\alpha c(\phi) \sum_{E \in E_h } \intE{ (\vr^k_h) ^{2 - \gamma} } \aleq h^\alpha.
\]
We may infer that
\begin{equation} \label{C7}
\frac{h^\alpha}{2} \left\| \sum_{\Gamma \in \Ghi} \intG{ \ju{\vr_h} \ju{\avo{\phi} } \chi \left( \frac{ \avg{\vu_h \cdot \vc{n} }}{h^\alpha} \right) }
\right\|_{L^1(0,T)} = \mathcal{O}(h^\beta), \beta > 0 \ \mbox{whenever}\ \alpha > 0.
\end{equation}

Next, using (\ref{S3}) again, we deduce
\[
\begin{split}
&\left| \sum_{E \in E_h} \sum_{\Gamma_E \subset \partial E} \int_{\Gamma_E} (\avo{\phi} - \phi) \ju{\vr^k_h} [ \avg{ \vu^k_h \cdot \vc{n} }]^- \ {\rm dS}_x
\right| \\ &\aleq h \| \phi \|_{C^1(\Ov{\Omega}_h)} \sum_{E \in E_h}\sum_{\Gamma_E \subset \partial E} \int_{\Gamma_E} | \ju{ (\vr^k_h)^{\gamma/2} } |
\ | \av{ (\vr^k_h)^{1 - \gamma/2} } | \ |\avg{ \vu^k_h \cdot \vc{n} }   | \ {\rm dS}_x\\
&\aleq h \left( \sum_{E \in E_h}\sum_{\Gamma_E \subset \partial E} \int_{\Gamma_E} \ju{ (\vr^k_h)^{\gamma/2} }^2 |\avg{ \vu^k_h \cdot \vc{n} }   | \ {\rm dS}_x
\right)^{1/2}\left( \sum_{E \in E_h}\sum_{\Gamma_E \subset \partial E} \int_{\Gamma_E} (\vr^k_h)^{2 - \gamma} |\vu^k_h| \ {\rm dS}_x
\right)^{1/2}\\
&\aleq h^{1/2} \left( \sum_{E \in E_h}\sum_{\Gamma_E \subset \partial E} \int_{\Gamma_E} \ju{ (\vr^k_h)^{\gamma/2} }^2 |\avg{ \vu^k_h \cdot \vc{n} }
 | \ {\rm dS}_x
\right)^{1/2}\left( \sum_{E \in E_h} \int_{E} (\vr^k_h)^{2 - \gamma} |\avo{\vu_h}| \ \dx;
\right)^{1/2}
\end{split}
\]
whence, using (\ref{S6}) to control the last term, we conclude
\begin{equation} \label{C8}
\left\| \sum_{E \in E_h} \sum_{\Gamma_E \subset \partial E} \int_{\Gamma_E} (\avo{\phi} - \phi) \ju{\vr_h} [ \avg{ \vu_h \cdot \vc{n} }]^- \ {\rm dS}_x
\right\|_{L^2(0,T)} = \mathcal{O}(h^\beta).
\end{equation}

Furthermore,
\[
\begin{split}
\sum_{E \in E_h} \sum_{\Gamma_E \subset \partial E} \int_{\Gamma_E} \phi \vr^k_h \Big(\vu^k_h \cdot \vc{n} - \avg{ \vu^k_h \cdot \vc{n} } \Big) \ {\rm dS}_x
= \sum_{E \in E_h} \sum_{\Gamma_E \subset \partial E} \int_{\Gamma_E} \left( \phi - \avg{\phi}
 \right) \vr^k_h \Big(\vu^k_h \cdot \vc{n} - \avg{ \vu^k_h \cdot \vc{n} } \Big) \ {\rm dS}_x,
\end{split}
\]
where, by virtue of Poincar\' e's inequality and the trace estimates (\ref{S4a}),
\[
\begin{split}
&\left|
\sum_{E \in E_h} \sum_{\Gamma_E \subset \partial E} \int_{\Gamma_E} \left( \phi - \avg{\phi} \right)
\vr^k_h \Big(\vu^k_h \cdot \vc{n} - \avg{ \vu^k_h \cdot \vc{n} } \Big) \ {\rm dS}_x \right| \\
&\aleq h \| \Grad \phi \|_{L^\infty(\Omega_h)}
\sum_{E \in E_h} \sum_{\Gamma_E \subset \partial E} \int_{\Gamma_E}  \vr^k_h \left| \vu^k_h \cdot \vc{n} - \avg{ \vu^k_h \cdot \vc{n} } \right| \ {\rm dS}_x
\aleq \sum_{E \in E_h}  \int_{E}  \vr^k_h \left| \vu^k_h \cdot \vc{n} - \avg{ \vu^k_h \cdot \vc{n} } \right| \ \dx\\
&\aleq h \sum_{E \in E_h} \| \nabla_h \vu^k_h \|_{L^2(E)} \| \vr^k_h \|_{L^2(E)} \aleq h \| \nabla_h \vu^k_h \|_{L^2(\Omega_h)} \| \vr^k_h \|_{L^2(\Omega_h)}
\aleq h \| \nabla_h \vu^k_h \|_{L^2(\Omega_h)} \| \vr^k_h \|_{L^\infty(\Omega_h)}^{1/2}.
\end{split}
\]
Going back to (\ref{S6}) we observe that the right-hand side is controlled as soon as
\begin{equation} \label{C9}
1 - \frac{2 + \alpha}{2 \gamma} > 0 \ \mbox{meaning}\  \alpha < 2 (\gamma - 1).
\end{equation}

Finally, it is easy to check that the last integral in (\ref{C6}) can be handled in the same way.
Thus we conclude that the consistency formulation (\ref{C10}) holds
for any test function $\varphi \in \DC([0, \infty) \times \Ov{\Omega}_h )$ as long as $\alpha > 0$, $\gamma >1$ are interrelated through (\ref{C9}).

\subsection{Consistency formulation of the momentum method}

Our goal is to take $\Pi^V_h [\vcg{\phi}]$, $\vcg{\phi} \in \DC(\Omega_h; R^3)$ as a test function in the momentum scheme (\ref{N4}). To begin, observe that
\[
\begin{split}
\intOh{ \nabla_h \vu_h : \nabla_h \Pi^V_h [\vcg{\phi}] } = \intOh{ \nabla_h \vu_h : \Grad \vcg{\phi}  },& \
\intOh{ {\rm div}_h \vu_h  {\rm div}_h \Pi^V_h [\vcg{\phi}] } = \intOh{ {\rm div}_h \vu_h  \Div \vcg{\phi} }\\
\intOh{ p(\vr_h)  {\rm div}_h \Pi^V_h [\vcg{\phi}] } &= \intOh{ p(\vr_h)  \Div \vcg{\phi} },
\end{split}
\]
see \cite[Chapter 9, Lemma 8]{FeiKaPok}.

\subsubsection{Time derivative}
\label{Tder}

We compute
\begin{equation} \label{C11}
\begin{split}
\intOh{ D_t (\vr^k_h \avo{\vu^k_h} ) \cdot \vcg{\phi} } &= \intOh{ D_t (\vr^k_h \avo{\vu^k_h}) \cdot \Pi^V_h [\vcg{\phi} ] }
\\
&+
\intOh{  \vr^{k-1}_h  \frac{ \avo{\vu^k_h} - \avo{\vu^{k-1}_h} }{\Delta t} \cdot \left( \vcg{\phi} - \Pi^V_h [\vcg{\phi}] \right) }
\\
&+ \intOh{ \frac{ \vr^k_h - \vr^{k-1}_h }{\Delta t} \avo{\vu^{k}_h} \cdot \left( \vcg{\phi} - \Pi^V_h [\vcg{\phi}] \right) },
\end{split}
\end{equation}
where
\[
\begin{split}
\left| \intOh{  \vr^{k-1}_h  \frac{ \avo{\vu^k_h} - \avo{\vu^{k-1}_h} }{\Delta t} \cdot \left( \vcg{\phi} - \Pi^V_h [\vcg{\phi}] \right) } \right|
&\aleq h^2 \| \phi \|_{C^2(\Ov{\Omega}_h)} \intOh{  \vr^{k-1}_h  \left| \frac{ \avo{\vu^k_h} - \avo{\vu^{k-1}_h} }{\Delta t} \right| }\\
\aleq h^2 \left( \intOh{ \vr^{k-1}_h } \right)^{1/2} & \left( \intOh{ \vr^{k-1}_h  \left( \frac{ \avo{\vu^k_h} - \avo{\vu^{k-1}_h} }{\Delta t} \right)^2  }
 \right)^{1/2}\\ \aleq h^2 (\TS)^{-1/2} &\left( \TS \intOh{ \vr^{k-1}_h  \left( \frac{ \avo{\vu^k_h} - \avo{\vu^{k-1}_h} }{\Delta t} \right)^2  }
 \right)^{1/2},
\end{split}
\]
where the most right integral is controlled in $L^2(0,T)$ by the numerical dissipation in (\ref{S3}).

As for the remaining integral, we may use inequality (\ref{C6a}) to obtain
\[
\begin{split}
&\left| \intOh{ \frac{ \vr^k_h - \vr^{k-1}_h }{\Delta t} \avo{\vu^{k}_h} \cdot \left( \vcg{\phi} - \Pi^V_h [\vcg{\phi}] \right) } \right|
\aleq h^2 \intOh{ \frac{ |\vr^k_h - \vr^{k-1}_h| }{\Delta t} | \avo{\vu^{k}_h} | }\\
&\aleq h^2 (\TS)^{-1} \| \vu^k_h \|_{L^6(\Omega_h; R^3)} \sup_{k} \| \vr^k_h \|_{L^{6/5}(\Omega_h)}
\aleq h^2 (\TS)^{-1} h^{-1/2} \| \vu^k_h \|_{L^6(\Omega_h; R^3)} \sup_{k} \| \vr^k_h \|_{L^{1}(\Omega_h)}
\end{split}
\]

Finally, we may repeat the same argument as in Section \ref{TD1} to conclude that
\begin{equation} \label{C13}
\begin{split}
\int_0^T &\intOh{ \psi D_t (\vr_h \avo{\vu_h} )  \Pi^V_h [\vcg{\phi} ] }\ \dt \\ &= - \int_0^T \intOh{ \vr_h \avo{\vu_h} \cdot \vcg{\phi} \partial_t \psi }
\ \dt
- \intOh{ \psi(0) \vr^0_h \avo{ \vu^0_h } \cdot \phi } + \mathcal{O}(h^\beta)
\end{split}
\end{equation}
provided $\psi \in \DC[0,T)$, $\vcg{\phi} \in \DC(\Omega_h; R^3)$.

\subsubsection{Convective term - upwind}

Applying formula (\ref{C4}) we obtain

\begin{equation} \label{C14}
\begin{split}
\intOh{ \vr^k_h \left( \avo{\vu^k_h} \otimes \vu^k_h \right) :  \Grad \vcg{\phi} }
&- \sum_{\Gamma \in \Ghi} \intG{ {\rm Up}[\vr^k_h \avo{\vu^k_h}, \vu^k_h] \cdot \ju{\avo{ \Pi^V_h[\vcg{\phi}] }} }
\\
&= \frac{h^\alpha}{2} \sum_{\Gamma \in \Ghi} \intG{ \ju{\vr^k_h \avo{ \vu^k_h } } \cdot \ju{ \avo{ \Pi^V_h[\vcg{\phi}] }}
\chi \left( \frac{ \avg{\vu^k_h \cdot \vc{n} }}{h^\alpha} \right) }
\\
&+ \sum_{E \in E_h} \sum_{\Gamma_E \subset \partial E} \int_{\Gamma_E} \left( \avo{ \Pi^V_h[\vcg{\phi}] } - \vcg{\phi} \right) \cdot
\ju{\vr^k_h \avo{ \vu^k_h }} [ \avg{ \vu^k_h \cdot \vc{n} }]^- \ {\rm dS}_x
\\
&+ \sum_{E \in E_h} \sum_{\Gamma_E \subset \partial E} \int_{\Gamma_E} \vr^k_h \vcg{\phi} \cdot \avo{ \vu^k_h } \Big(\vu^k_h \cdot \vc{n} - \avg{ \vu^k_h \cdot \vc{n} } \Big) \ {\rm dS}_x
\\
&+ \intOh{ \vr^k_h \avo{ \vu^k_h } \cdot \left( \avo{ \Pi^V_h[\vcg{\phi}] } - \vcg{\phi} \right) \Divh \vu^k_h }.
\end{split}
\end{equation}

We proceed in several steps.

\medskip

{\bf Step 1}

\medskip

Applying (\ref{S7}) we get
\[
\begin{split}
&\left| \frac{h^\alpha}{2} \sum_{\Gamma \in \Ghi} \intG{ \ju{\vr^k_h \avo{ \vu^k_h } } \cdot \ju{ \avo{ \Pi^V_h[\vcg{\phi}] }}
\chi \left( \frac{ \avg{\vu^k_h \cdot \vc{n} }}{h^\alpha} \right) } \right|\\
&\aleq h^{1 + \alpha}  \sum_{\Gamma \in \Ghi}
\intG{ \left|\ \ju{\vr^k_h \avo{ \vu^k_h } } \right|\ \chi \left( \frac{ \avg{\vu^k_h \cdot \vc{n} }}{h^\alpha} \right)} ,
\end{split}
\]
where
\begin{equation} \label{C15}
\ju{\vr^k_h \avo{ \vu^k_h } } = (\vr^k_h)^{\rm out} \ju {\avo{ \vu^k_h }} + \avo{\vu^k_h} \ju{ \vr^k_h }.
\end{equation}
Consequently
\[
\begin{split}
&\left| \frac{h^\alpha}{2} \sum_{\Gamma \in \Ghi} \intG{ \ju{\vr^k_h \avo{ \vu^k_h } } \cdot \ju{ \avo{ \Pi^V_h[\vcg{\phi}] }}
\chi \left( \frac{ \avg{\vu^k_h \cdot \vc{n} }}{h^\alpha} \right) } \right| \\
&\aleq h^{1 + \alpha} \left( \sum_{\Gamma \in \Ghi}
\intG{ \av{\vr^k_h} \ju{ \avo{ \vu^k_h } }^2 \chi \left( \frac{ \avg{\vu^k_h \cdot \vc{n} }}{h^\alpha} \right)} \right)^{1/2} \left( \sum_{\Gamma \in \Ghi}
\intG{ \vr^k_h  } \right)^{1/2} \\
&+ h^{1 + \alpha}  \sum_{\Gamma \in \Ghi}
\intG{ \left| \avo{\vu^k_h}  \ju{ \vr^k_h } \right| } ,
\end{split}
\]
where the first integral on the right-hand side is controlled by the numerical dissipation in (\ref{S3}) and the trace estimates.

Finally, applying the inequality (\ref{S4}), trace inequality (\ref{S4a})  and Sobolev's inequality (\ref{S6b}), we obtain
\[
\begin{split}
h^{1 + \alpha} & \sum_{\Gamma \in \Ghi}
\intG{ \left|\ \avo{\vu^k_h}  \ju{ \vr^k_h }  \ \right| }  \aleq h^{1 + \alpha}\sum_{\Gamma \in \Ghi}
\intG{ |\avo{\vu^k_h}|\ \left| \av{ \vr^k_h }^{1 - \gamma/2} \right| \   \left| \ju{ (\vr^k_h)^{\gamma/2}  } \right|  } \\
&\aleq h^{1 + \alpha}  \sum_{\Gamma \in \Ghi} \left( \intG{ \ju{ (\vr^k_h)^{\gamma/2}  }^2 } \right)^{1/2}
\| \avo{ \vu^k_h } \|_{L^6(\Gamma)} \left\| (\vr^k_h)^{1 - \gamma/2} \right\|_{L^3 (\Gamma)} \\
&\aleq h^{\frac{1 + \alpha}{2}}  \sum_{E \in E_h} \left( h^\alpha \int_{\partial E} \ju{ (\vr^k_h)^{\gamma/2}  }^2 \ {\rm dS}_x \right)^{1/2}
\| \avo{ \vu^k_h } \|_{L^6(E)} \left\| (\vr^k_h)^{1 - \gamma/2} \right\|_{L^3 (E)}\\
&\aleq h^{\frac{1 + \alpha}{2}} \left\| \nabla_h \vu_h \right\|_{L^2(\Omega_h)} \left\| (\vr^k_h)^{1 - \gamma/2} \right\|_{L^3 (\Omega_h)},
\end{split}
\]
where we have used the numerical dissipation in (\ref{S3}). Thus, in order to complete the estimates we have to control
\[
\left\| (\vr^k_h)^{1 - \gamma/2} \right\|_{L^3 (\Omega_h)}
\]
uniformly in $k$. As $1 < \gamma < 2$, it is enough to consider the critical case $\gamma = 1$, for which the inverse inequality
(\ref{S4})  gives rise to
\[
\left\| (\vr^k_h)^{1/2} \right\|_{L^3 (\Omega_h)} = \left( \left\| (\vr^k_h) \right\|_{L^{3/2 (\Omega_h)}} \right)^{1/2}
\aleq h^{-1/2} \| \vr^k_h \|^{1/2}_{L^1(\Omega_h)}.
\]

\medskip

{\bf Step 2}

\medskip

Using (\ref{C15}) we deduce
\[
\begin{split}
\sum_{E \in E_h} \sum_{\Gamma_E \subset \partial E} \int_{\Gamma_E} & \left( \avo{ \Pi^V_h[\vcg{\phi}] } - \vcg{\phi} \right) \cdot
\ju{\vr^k_h \avo{ \vu^k_h }} [ \avg{ \vu^k_h \cdot \vc{n} }]^- \ {\rm dS}_x\\
&= \sum_{E \in E_h} \sum_{\Gamma_E \subset \partial E} \int_{\Gamma_E}  \left( \avo{ \Pi^V_h[\vcg{\phi}] } - \vcg{\phi} \right) \cdot
\Big( (\vr^k_h)^{\rm out} \ju {\avo{ \vu^k_h }} + \avo{\vu^k_h} \ju{ \vr^k_h } \Big) [ \avg{ \vu^k_h \cdot \vc{n} }]^- \ {\rm dS}_x
\end{split},
\]
where, furthermore,
\[
\begin{split}
&\left| \sum_{E \in E_h} \sum_{\Gamma_E \subset \partial E} \int_{\Gamma_E} \left( \avo{ \Pi^V_h[\vcg{\phi}] } - \vcg{\phi} \right)
(\vr^k_h)^{\rm out} \ju {\avo{ \vu^k_h } } [ \avg{ \vu^k_h \cdot \vc{n} }]^- \ {\rm dS}_x \right|\\ &\aleq
h^2 \| \vcg{\phi} \|_{C^2(\Ov{\Omega}_h;R^3)} \left( \sum_{E \in E_h} \sum_{\Gamma_E \subset \partial E} \int_{\Gamma_E} -
(\vr^k_h)^{\rm out} \ju {\avo{ \vu^k_h } }^2 [ \avg{ \vu^k_h \cdot \vc{n} }]^- \ {\rm dS}_x
\right)^{1/2} \times \\
& \times \left( \sum_{E \in E_h} \sum_{\Gamma_E \subset \partial E} \int_{\Gamma_E}
(\vr^k_h)^{\rm out} |\avo{\vu^k_h}| \ {\rm dS}_x \right)^{1/2},
\end{split}
\]
where the former integral in the product on the right-hand is controlled by the numerical dissipation in (\ref{S3}), while
\[
\sum_{E \in E_h} \sum_{\Gamma_E \subset \partial E} \int_{\Gamma_E}
(\vr^k_h)^{\rm out} |\avo{\vu^k_h}| \ {\rm dS}_x \aleq h^{-1} \| \vu^k_h \|_{L^6(\Omega_h; R^3)} \| \vr^k_h \|_{L^{6/5} (\Omega_h)},
\aleq h^{-3/2} \| \vu^k_h \|_{L^6(\Omega_h; R^3)} \| \vr^k_h \|_{L^{1} (\Omega_h)}.
\]

Finally,
\[
\begin{split}
&\left| \sum_{E \in E_h} \sum_{\Gamma_E \subset \partial E} \int_{\Gamma_E}  \left( \avo{ \Pi^V_h[\vcg{\phi}] } - \vcg{\phi} \right) \cdot
 \avo{\vu^k_h} \ju{ \vr^k_h }  [ \avg{ \vu^k_h \cdot \vc{n} }]^- \ {\rm dS}_x \right| \\ &\aleq h^2
 \sum_{E \in E_h} \sum_{\Gamma_E \subset \partial E} \| \avo{\vu^{k}_h} \|_{L^6(\Gamma)}^2 \| \vr^k_h \|_{L^{3/2}(\Gamma)}
 \aleq h \| \vu^k_h \|^2_{L^6(\Omega_h)} \| \vr^k_h \|_{L^{3/2} (\Omega_h)} \aleq h^{3 - 3/\gamma} \| \vu^k_h \|^2_{L^6(\Omega_h)} \| \vr^k_h \|_{L^{\gamma}
 (\Omega_h)},
\end{split}
\]
where the exponent $3 - 3/\gamma > 0$ as soon as $\gamma > 1$.

\medskip

{\bf Step 3}

\medskip

We write
\[
\begin{split}
\sum_{E \in E_h} \sum_{\Gamma_E \subset \partial E} &\int_{\Gamma_E} \vr^k_h \vcg{\phi} \cdot \avo{ \vu^k_h }
\Big(\vu^k_h \cdot \vc{n} - \avg{ \vu^k_h \cdot \vc{n} } \Big) \ {\rm dS}_x \\ &=
\sum_{E \in E_h} \sum_{\Gamma_E \subset \partial E} \int_{\Gamma_E} \vr^k_h ( \vcg{\phi} - \avg{ \vcg{\phi} } )  \cdot \avo{ \vu^k_h }
\Big(\vu^k_h \cdot \vc{n} - \avg{ \vu^k_h \cdot \vc{n} } \Big) \ {\rm dS}_x,
\end{split}
\]
where, by virtue of the trace inequality (\ref{S4a}) and Poincar\' e's inequality (\ref{S4b}),
\[
\begin{split}
&\left| \sum_{E \in E_h} \sum_{\Gamma_E \subset \partial E} \int_{\Gamma_E} \vr^k_h ( \vcg{\phi} - \avg{ \vcg{\phi} } )  \cdot \avo{ \vu^k_h }
\Big(\vu^k_h \cdot \vc{n} - \avg{ \vu^k_h \cdot \vc{n} } \Big) \ {\rm dS}_x \right| \\
&\aleq h  \left\| \sqrt{\vr^k_h} \right\|_{L^\infty(\Omega_h)} \sum_{E \in E_h} \sum_{\Gamma_E \subset \partial E} \int_{\Gamma_E} \sqrt{\vr^k_h} | \avo{ \vu^k_h } |
\Big| \vu^k_h \cdot \vc{n} - \avg{ \vu^k_h \cdot \vc{n} } \Big| \ {\rm dS}_x\\
&\aleq h  \left\| \sqrt{\vr^k_h} \right\|_{L^\infty(\Omega_h)} \| \sqrt{\vr^{k}_h} \avo{ \vu^k_h } \|_{L^2(\Omega_h)} \|
\nabla_h \vu^k_h \|_{L^2(\Omega_h; R^3)},
\end{split}
\]
where, in view of (\ref{S6})
\[
\left\| \sqrt{\vr^k_h} \right\|_{L^\infty(\Omega_h)} \aleq h^{- \frac{2 + \alpha}{2 \gamma}},
\ \mbox{with} \frac{2 + \alpha}{2 \gamma} < 1 \ \mbox{or}\ 0 < \alpha < 2 (\gamma - 1).
\]

\medskip

{\bf Step 4}

\medskip

Finally,
\[
\begin{split}
&\left| \intOh{ \vr^k_h \avo{ \vu^k_h } \cdot \left( \avo{ \Pi^V_h[\vcg{\phi}] } - \vcg{\phi} \right) \Divh \vu^k_h } \right| \\ &\aleq
h^2  \left\| \sqrt{\vr^k_h} \right\|_{L^\infty(\Omega_h)}  \| \sqrt{\vr^{k}_h} \avo{ \vu^k_h } \|_{L^2(\Omega_h)} \|
\nabla_h \vu^k_h \|_{L^2(\Omega_h; R^3)};
\end{split}
\]
whence the rest of the proof follows exactly as in Step 3.

Summing up the previous observations, we obtain the consistency formulation of the momentum method (\ref{C16}).

\begin{Remark} \label{RC1}

As $\vcg{\varphi}$ has compact support, equation (\ref{C16}) is satisfied also on the limit domain $\Omega$ for all $h$ small enough.

\end{Remark}

Thus we have shown the following result.

\begin{Proposition} \label{PC1}

Let the pressure $p$ satisfy (\ref{i5}), with $1 < \gamma < 2$. Suppose that $[ \vr_h, \vu_h]$ is a family of numerical solutions
given through (\ref{C1}), (\ref{C2}), where $[\vr^k_h, \vu^k_h]$ satisfy (\ref{N2}--\ref{N4}), where
\begin{equation} \label{expo}
\TS \approx h,\ 0 < \alpha < 2 (\gamma - 1).
\end{equation}

Then
\[
- \intOh{ \vr^0_h \varphi (0, \cdot) } = \int_0^T \intOh{ \left[  \vr_h \partial_t \varphi + \vr_h \vu_h \cdot \Grad \varphi \right] } \dt
+ \mathcal{O}(h^\beta), \ \beta > 0,
\]
for any test function $\varphi \in \DC([0, \infty) \times \Ov{\Omega}_h )$,
\begin{equation} \label{E2}
\begin{split}
- \intOh{ \vr^0_h \avo{ \vu^0_h } \cdot  \vcg{\varphi}(0, \cdot) } &= \int_0^T \intOh{
\Big[ \vr_h \avo{\vu_h} \partial_t \vcg{\varphi} + \vr_h \avo{\vu_h} \otimes \vu_h : \Grad \vcg{\varphi} + p(\vr_h ) \Div \vcg{\varphi} \Big] } \ \dt\\
&- \int_0^T \intOh{ \Big[ \mu \nabla_h \vu_h : \Grad \vcg{\varphi} + (\mu/3 + \eta) {\rm div}_h \vu_h \cdot \Div \vcg{\varphi} \Big] } \dt
+ \mathcal{O}(h^\beta), \ \beta > 0,
\end{split}
\end{equation}
for any $\vcg{\varphi} \in \DC([0, \infty) \times \Omega_h; R^3)$.

Moreover, the solution satisfies the energy inequality
\begin{equation} \label{E3}
\begin{split}
\intOh{ \left[ \frac{1}{2} \vr_h | \avo{\vu_h } |^2 + P(\vr_h) \right] (\tau, \cdot)} &+
\int_0^\tau \intOh{ \mu |\nabla_h \vu_h |^2 + (\mu/3 + \eta) |{\rm div}_h \vu_h |^2 } \ \dt \\ &\leq
\intOh{ \left[ \frac{1}{2} \vr^0_h | \avo{\vu^0_h } |^2 + P(\vr^0_h) \right]}
\end{split}
\end{equation}
for a.e. $\tau \in [0,T]$.

\end{Proposition}

\begin{Remark} \label{RC3}

A close inspection of the previous discussion shows that the same method can be used to handle a variable time step $\TS_k$ adjusted for each step of iteration by means of a CFL-type condition, {\cred such as $ || \vu_h^{k-1} + c_h^{k-1} ||_{L^\infty(\Omega)} \TS_k / h \leq CFL $. Here $CFL \in (0,1]$ and  $c_h^{k-1} \equiv \sqrt{p'(\rho_h^{k-1})}$ denotes the sound speed. Though this condition is necessary for stability of time-explicit
numerical schemes, it still may be appropriate even for implicit schemes for areas of high-speed flows. Note
that the only part that must be changed in the {proof of Proposition \ref{PC1} is Section \ref{Tder}}, where the time derivative in the momentum method
is estimated.}
\end{Remark}

\section{Measure-valued solutions}
\label{M}

Our ultimate goal is to perform the limit $h \to 0$. For the sake of simplicity, we consider the initial data
\[
\vr_0 \in L^\infty( R^3), \ \vr^0 \geq \underline{\vr} > 0 \ \mbox{a.a. in}\ R^3,\
\vu_0 \in L^2 (R^3).
\]
With this ansatz, it is easy to find the approximation $[\vr^0_h, \vu^0_h]$ such that
\begin{equation} \label{M1}
\begin{split}
\vr^0_h &\to \vr_0 \ \mbox{in}\ L^\gamma_{\rm loc}(\Omega), \ \vr^0_h > 0, \ \intOh{\vr^0_h \phi} \to
\intO{ \vr_0 \phi } \ \mbox{for any}\ \phi \in L^\infty(R^3),\\
\vr^0_h \avo{\vu^0_h} &\to \vr_0 \vu_0 \ \mbox{in}\ L^2_{\rm loc}(\Omega; R^3),\
\intOh{ \vr^0_h \avo{ \vu^0_h } \cdot \phi } \to \intO{ \vr_0 \vu_0  \cdot \phi }
\ \mbox{for any} \ \phi \in L^\infty(R^3; R^3),\\
&\intOh{ \left[ \frac{1}{2} \vr^0_h |\avo{\vu^0_h} |^2 + P(\vr^0_h) \right] } \to
\intO{ \left[ \frac{1}{2} \vr_0 |{\vu_0} |^2 + P(\vr_0) \right] } \ \mbox{as}\ h \to 0.
\end{split}
\end{equation}

\subsection{Weak limit}

Extending $\vr_h$ by $\underline{\vr} > 0$ and $\vu_h$ to be zero outside $\Omega_h$, we may use the energy estimates (\ref{E3}) to
deduce that, at least for suitable subsequences,
\[
\begin{split}
\vr_h &\to \vr \ \mbox{weakly-(*) in}\ L^\infty(0,T; L^\gamma(\Omega)),\ \vr \geq 0 \\
\avo{\vu_h}, \ \vu_h &\to \vu \ \mbox{weakly in}\ L^2((0,T) \times \Omega; R^3),\\
\mbox{where} \ \vu &\in L^2(0,T; W^{1,2}_0(\Omega)), \ \nabla_h \vu_h \to \Grad \vu \ \mbox{weakly in} \ L^2((0,T) \times \Omega;
R^{3 \times 3}), \\
\vr_h \avo{\vu_h}  &\to \Ov{\vr_h \vu_h} \ \  \mbox{weakly-(*) in}\ L^\infty(0,T; L^{\frac{2\gamma}{\gamma + 1}}(\Omega; R^3)),
\end{split}
\]
see { \cite{FeKaMi} or \cite[Part II, Section 10.4]{FeiKaPok}}.

\begin{Remark} \label{RM1}

Note that, by virtue of Poincar\' e's inequality (\ref{S4b}) and the energy estimates (\ref{E3}),
\[
\| \vu_h - \avo{\vu_h } \|_{L^2(0,T; L^2(K; R^3)) } \aleq h \ \mbox{for any compact}\ K \in \Omega,
\]
in particular, the weak limits of $\vu_h$, $\avo{\vu_h}$ coincide in $\Omega$.

\end{Remark}

In addition, the limit functions satisfy the equation of continuity in the form
\begin{equation} \label{M2}
- \intO{ \vr_0 \varphi (0, \cdot) } = \int_0^T \intO{ \left[  \vr \partial_t \varphi + \Ov{\vr \vu} \cdot \Grad \varphi \right] } \dt
\end{equation}
for any test function $\varphi \in \DC([0, \infty) \times \Ov{\Omega} )$. It follows from (\ref{M2}) that
$\vr \in C_{\rm weak}([0,T]; L^\gamma(\Omega))$; whence (\ref{M2}) can be rewritten as
\begin{equation} \label{M3}
\left[ \intO{ \vr \varphi (\tau, \cdot) } \right]_{t = 0}^{t = \tau} = \int_0^\tau
\intO{ \left[  \vr \partial_t \varphi + \Ov{\vr \vu} \cdot \Grad \varphi \right] } \dt
\end{equation}
for any $0 \leq \tau \leq T$ and any $\varphi \in C^\infty([0,T] \times \Ov{\Omega})$.

\subsection{Young measure generated by numerical solutions}
\label{young}

The energy inequality (\ref{S1}), along with the consistency (\ref{C10}), (\ref{C16}) provide a suitable platform for the use of the theory of
measure-valued solutions developed in \cite{FGSWW1}. Consider the family $[\vr_h, \vu_h]$. In accordance with the weak convergence statement
derived in the preceding part, this family generates a Young measure - a parameterized measure
\[
\nu_{t,x} \in L^\infty((0,T) \times \Omega; \mathcal{P}([0, \infty) \times R^3)) \ \mbox{for a.a.}\ (t,x) \in (0,T) \times \Omega,
\]
such that
\[
\left< \nu_{t,x}, g(\vr, \vu) \right> = \Ov{g(\vr, \vu)}(t,x)\ \mbox{for a.a.}\ (t,x) \in (0,T) \times \Omega,
\]
whenever $g \in C([0, \infty) \times R^3)$, and
\[
g(\vr_h, \vu_h) \to \Ov{g(\vr, \vu)} \ \mbox{weakly in}\ L^1((0,T) \times \Omega).
\]
Moreover, in view of Remark \ref{RM1}, the Young measures generated by $[\vr_h, \vu_h]$ and $[\vr, \avo{\vu_h}]$ coincide for a.a.
$(t,x) \in (0,T) \times \Omega$.

Accordingly, the equation of continuity (\ref{M3}) can be written as
\begin{equation} \label{M4}
\left[ \intO{ \vr \varphi (\tau, \cdot) } \right]_{t = 0}^{t = \tau} = \int_0^\tau
\intO{ \left[  \vr \partial_t \varphi + \left< \nu_{t,x}, \vr \vu \right>  \cdot \Grad \varphi \right] } \dt
\end{equation}

In order to apply a similar treatment to the momentum equation (\ref{E2}), { we have to replace the expression
$\vr_h \avo{\vu_h} \otimes \vu_h$ in the convective term by $\vr_h \avo{\vu_h} \otimes \avo{\vu_h}$.} This is possible as
\[
\begin{split}
& \left\| \vr_h \avo{\vu_h} \otimes \vu_h - \vr_h \avo{\vu_h} \otimes \avo{\vu_h} \right\|_{L^1(\Omega_h; R^{3 \times 3})}
=\left\| \vr_h \avo{\vu_h} \otimes (\vu_h -  \avo{\vu_h}) \right\|_{L^1(\Omega_h; R^{3 \times 3})}\\
&\aleq h \| \sqrt{\vr_h} \avo{\vu_h} \|_{L^2(\Omega_h); R^3)} \| \nabla_h \vu_h \|_{L^2(\Omega_h; R^{3 \times 3})}
\| \sqrt{\vr_h} \|_{L^\infty(\Omega_h)},
\end{split}
\]
where, by virtue of (\ref{S6}),
\[
h \| \sqrt{\vr_h} \|_{L^\infty(\Omega_h)} \aleq h^{ 1 - \frac{2 + \alpha}{2 \gamma}},
\]
where the exponent is positive { as soon as (\ref{expo}) holds, specifically, $0 < \alpha < 2 (\gamma - 1)$}. Moreover, we have
\[
\vr_h \avo{\vu_h} \otimes \avo{ \vu_h} + p(\vr_h) \mathbb{I} \to
\left\{ \vr \vu \otimes \vu + p(\vr) \mathbb{I} \right\}\ \mbox{weakly-(*) in}\ \left[ L^\infty(0,T; \mathcal{M}(\Omega) )\right]^{3 \times 3};
\]
whence letting $h \to 0$ in (\ref{E2}) gives rise to
\[
\begin{split}
- \intO{ \vr_0 \vu_0 \cdot  \vcg{\varphi}(0, \cdot) } &= \int_0^T \intO{
\Big[ \left< \nu_{t,x}; \vr \vu \right> \partial_t \vcg{\varphi} + \left\{ \vr \vu \otimes \vu + p(\vr) \mathbb{I} \right\} : \Grad
\vcg{\varphi}  \Big] } \ \dt\\
&- \int_0^T \intOh{ \Big[ \mu \nabla \vu : \Grad \vcg{\varphi} + (\mu/3 + \eta) {\rm div} \vu \cdot \Div \vcg{\varphi} \Big] } \dt
\end{split}
\]
or, equivalently,
\begin{equation} \label{M5}
\begin{split}
\left[ \intO{ \left< \nu_{t,x} ; \vr \vu \right> \cdot  \vcg{\varphi}(0, \cdot) } \right]_{t = 0}^{t = \tau} &= \int_0^\tau \intO{
\Big[ \left< \nu_{t,x}; \vr \vu \right> \cdot \partial_t \vcg{\varphi} + \left\{ \vr \vu \otimes \vu + p(\vr) \mathbb{I} \right\} : \Grad
\vcg{\varphi}  \Big] } \ \dt\\
&- \int_0^\tau \intO{ \Big[ \mu \nabla \vu : \Grad \vcg{\varphi} + (\mu/3 + \eta) {\rm div} \vu \cdot \Div \vcg{\varphi} \Big] } \dt
\end{split}
\end{equation}
for any $0 \leq \tau \leq T$, $\varphi \in \DC([0,T] \times \Omega; R^3)$, where we have set
\[
\nu_{0,x} = \delta_{[\vr_0(x), \vu_0(x)]}.
\]

Finally, we introduce the \emph{concentration remainder}
\[
\mathcal{R} = \left\{ \vr \vu \otimes \vu + p(\vr) \mathbb{I} \right\} - \left< \nu_{t,x}; \vr \vu \otimes \vu + p(\vr) \mathbb{I} \right>
\in [ L^\infty(0,T; \mathcal{M}(\Omega)) ]^{3 \times 3}
\]
and rewrite (\ref{M5}) in the form
\begin{equation} \label{M6}
\begin{split}
&\left[ \intO{ \left< \nu_{t,x} ; \vr \vu \right> \cdot  \vcg{\varphi}(0, \cdot) } \right]_{t = 0}^{t = \tau} \\ &= \int_0^\tau \intO{
\Big[ \left< \nu_{t,x}; \vr \vu \right> \cdot \partial_t \vcg{\varphi} + \left< \nu_{t,x}; \vr \vu \otimes \vu \right>
: \Grad \vcg{\varphi} + \left< \nu_{t,x}, p(\vr) \right> \Div
\vcg{\varphi}  \Big] } \ \dt\\
&- \int_0^\tau \intO{ \Big[ \mu \nabla \vu : \Grad \vcg{\varphi} + (\mu/3 + \eta) {\rm div} \vu \cdot \Div \vcg{\varphi} \Big] } \dt
+ \int_0^\tau \intO{ \mathcal{R} : \Grad \varphi } \dt
\end{split}
\end{equation}
for any $0 \leq \tau \leq T$, $\varphi \in \DC([0,T] \times \Omega; R^3)$.

Similarly, the energy inequality (\ref{E3}) can be written as
\begin{equation} \label{M7}
\begin{split}
\left[ \intO{ \left[ \frac{1}{2} \left< \nu_{t,x}; \vr | {\vu } |^2 + P(\vr) \right> \right] } \right]_{t = 0}^{t = \tau} &+
\int_0^\tau \intOh{ \mu |\nabla \vu |^2 + (\mu/3 + \eta) |{\rm div} \vu |^2 } \ \dt \\ + \mathcal{D}(\tau) &\leq 0
\end{split}
\end{equation}
for a.e. $\tau \in [0,T]$, with the \emph{dissipation defect} $\mathcal{D}$ satisfying
\begin{equation} \label{M8}
\int_0^\tau \| \mathcal{R} \|_{\mathcal{M}(\Omega)}\ \dt \aleq \int_0^\tau \mathcal{D}(t) \ \dt,\
\mathcal{D}(\tau) \geq  \liminf_{h \to \infty} \int_0^\tau \intOh{ |\nabla_h \vu_h|^2 }\ \dt - \int_0^\tau \intO{ |\Grad \vu |^2 }\ \dt,
\end{equation}
cf. \cite[Lemma 2.1]{FGSWW1}.

{At this stage, we recall the concept of \emph{dissipative measure valued solution} introduced in \cite{FGSWW1}. These are measure--valued
solutions of the Navier-Stokes system (\ref{i1}--\ref{i4}) satisfying the energy inequality (\ref{M7}), where the concentration remainder
in the momentum equation is dominated by the dissipation defect as stated in (\ref{M8}) and the following analogue of Poincar\' e's inequality
holds:
\begin{equation} \label{PI}
\lim_{h \to 0} \int_0^\tau \intOh{ | \vu_h - \vu |^2 } \ \dt \leq \liminf_{h \to \infty} \int_0^\tau \intOh{ |\nabla_h \vu_h|^2 }
\dt - \int_0^\tau \intO{ |\Grad \vu |^2 } \ \dt
(\leq \mathcal{D}(\tau)) ,
\end{equation}
where $\vu$ is a weak limit of $\vu_h$, or, equivalently, of $\avo{\vu_h}$.
Consequently,
relations (\ref{M4}), (\ref{M6}--\ref{M8}) imply that the Young measure $\{ \nu_{t,x} \}_{t,x \in (0,T) \times \Omega}$ represents
a dissipative measure-valued solution of the Navier-Stokes system (\ref{i1}--\ref{i4})
in the sense of \cite{FGSWW1} as soon as we check (\ref{PI}).}

By standard Poincar\' e's inequality in $\Omega_h$ we get, on one hand,
\[
\intOh{ |\vu_h - \vu |^2 } = \intOh{ |\vu_h - \Pi^V_h [\vu] |^2 } + \intOh{ |\Pi^V_h [\vu] - \vu |^2 }
\aleq \intOh{ |\nabla_h \vu_h - \nabla_h \Pi^V_h \vu |^2 } + \mathcal{O}(h^\beta).
\]
On the other hand,
\[
\liminf_{h \to \infty} \int_0^\tau \intOh{ |\nabla_h \vu_h|^2 } \dt -  \int_0^\tau \intO{ |\Grad \vu |^2 }\ \dt =
\liminf_{h \to \infty} \int_0^\tau \intOh{ |\nabla_h \vu_h - \Grad \vu |^2 } \dt.
\]
Thus it is enough to observe that, by virtue of (\ref{S7}),
\[
\nabla_h \Pi^V_h [\vu] \to \Grad \vu \ \mbox{(strongly) in}\ L^2(\Omega_h;R^3) \ \mbox{whenever}\ \vu \in W^{1,2}_0 (\Omega; R^3).
\]

Seeing that validity of (\ref{M6}) as well as the bound on the dissipation remainder (\ref{M8}) can be extended to the class of test functions
$\varphi \in C^1([0,T] \times \Ov{\Omega}; R^3)$, $\varphi|_{\partial \Omega} = 0$,
we have shown the following result.

\begin{Theorem} \label{MT1}
Let the pressure $p$ satisfy (\ref{i5}), with $1 < \gamma < 2$. Suppose that $[ \vr_h, \vu_h]$ is a family of numerical solutions
given through (\ref{C1}), (\ref{C2}), where $[\vr^k_h, \vu^k_h]$ satisfy (\ref{N2}--\ref{N4}), where
\[
\TS \approx h,\ 0 < \alpha < 2 (\gamma - 1),
\]
and the initial data satisfy (\ref{M1}).

Then any Young measure $\{ \nu_{t,x} \}_{t,x \in (0,T) \times \Omega}$ generated by  $[\vr^k_h, \vu^k_h]$ for $h \to 0$ represents a
dissipative measure-valued solution of the Navier-Stokes system (\ref{i1}--\ref{i4}) in the sense of \cite{FGSWW1}.

\end{Theorem}

Of course, the conclusion of Theorem \ref{MT1} is rather weak, and, in addition, the Young measure need not be unique. On the other hand, however,
we may use the weak-strong uniqueness principle established in \cite[Theorem 4.1]{FGSWW1} to obtain our final convergence result.
\begin{Theorem} \label{MT2}
In addition to the hypotheses of Theorem \ref{MT1}, suppose that the Navier-Stokes system (\ref{i1}--\ref{i4}) endowed with the initial
data $[\vr_0,  \vu_0]$ admits a regular solution $[\vr, \vu]$ belonging to the class
\[
\vr, \ \Grad \vr,\ \vu, \Grad \vu \in C([0,T] \times \Ov{\Omega}), \ \partial_t \vu \in L^2(0,T; C(\Ov{\Omega}; R^3)),\
\vr > 0, \ \vu|_{\partial \Omega} = 0.
\]

Then
\[
\vr_h \to \vr \ \mbox{(strongly) in}\ L^\gamma((0,T) \times K), \ \vu_h \to \vu \ \mbox{(strongly) in}\ L^2((0,T) \times K; R^3)
\]
for any compact $K \subset \Omega$.

\end{Theorem}

{Indeed, the weak--strong uniqueness implies that the Young measure generated by the family of numerical solutions coincides
at each point $(t,x)$ with the Dirac mass supported by the smooth solution of the problem. In particular, the numerical solutions converge strongly and
no oscillations occur. Note that the Navier--Stokes system admits local-in-time strong solutions for arbitrary smooth initial data, see
e.g. Cho et al. \cite{ChoChoeKim} , and even global-in-time smooth solutions for small initial data, see, e.g.,
Matsumura and Nishida \cite{MANI}, as soon as the physical domain $\Omega$ is sufficiently smooth. }

{\cred
\section{Conclusions}
We have studied the convergence of numerical solutions obtained by
the mixed finite element--finite volume scheme applied to the isentropic Navier-Stokes equations.
We have assumed the isentropic pressure--density
state equation $p(\vr)=a \vr^\gamma$ with $\gamma \in (1,2)$. Remind that
this assumption is not restrictive, since the largest physically relevant
exponent is $\gamma = 5/3$. In order to establish the convergence result we have used
the concept of dissipative measure-valued solutions. These are the measure-valued solutions, that, in addition, satisfy
an energy inequality in which the dissipation defect measure
dominates the concentration remainder in the equations.
The energy inequality (\ref{S1}), along with the consistency (\ref{C10}), (\ref{C16}) gave us a suitable framework to apply the theory of
measure-valued solutions. As shown in Section~\ref{young} the numerical solutions $[\vr_h, \vu_h]$ generate a Young measure - a parameterized measure $\{ \nu_{t,x} \}_{t,x \in (0,T) \times \Omega}$, that represents a dissipative measure-valued solution of the Navier-Stokes system (\ref{i1}--\ref{i4}), cf.~Theorem \ref{MT1}. Finally, using
the weak-strong uniqueness principle established in \cite[Theorem 4.1]{FGSWW1} we have obtained the convergence of the numerical solutions to the exact regular solution, as long as the latter exists, cf.~Theorem~\ref{MT2}.
The present result is the first convergence result for numerical solutions of three-dimensional compressible isentropic
Navier-Stokes equations in the case of full adiabatic exponent $\gamma \in (1,2).$
}

\def\cprime{$'$} \def\ocirc#1{\ifmmode\setbox0=\hbox{$#1$}\dimen0=\ht0
  \advance\dimen0 by1pt\rlap{\hbox to\wd0{\hss\raise\dimen0
  \hbox{\hskip.2em$\scriptscriptstyle\circ$}\hss}}#1\else {\accent"17 #1}\fi}


\end{document}